\documentclass[11pt]{article}
\usepackage{indentfirst}
\usepackage{dsfont}
\usepackage{amssymb}
\usepackage{epsfig}

\def\be{\begin{equation}}
\def\ee{\end{equation}}
\def\bea{\begin{eqnarray}}
\def\eea{\end{eqnarray}}
\def\bes{\begin{eqnarray*}}
\def\ees{\end{eqnarray*}}

\def\nn{\nonumber}
\def\<{\langle}
\def\>{\rangle}
\def\lb{\label}
\def\bs{\setminus}

\def\R{{\bf R}}
\def\C{{\bf C}}
\def\Z{{\bf Z}}

\def\U{{\bf U}}

\def\ga{{\gamma}}

\def\th{{\theta}}
\def\om{{\omega}}

\def\ep{{\epsilon}}

\def\sg{{\sigma}}

\def\Sg{{\Sigma}}

\def\H{{\cal H}}

\def\P{{\cal P}}
\def\J{{\cal J}}

\def\Sp{{\rm Sp}}

\def\dm{{\rm \diamond}}

\def\hb{\vrule height0.18cm width0.14cm $\,$}

\title{ Elliptic and non-hyperbolic closed characteristics on compact convex P-cyclic symmetric hypersurfaces in ${\bf R}^{2n}$}
\author{Hui Liu$^{1}$, \thanks{Partially supported by NSFC (Nos. 11401555, 11771341).
E-mail: huiliu00031514@whu.edu.cn}\qquad  Chongzhi Wang$^{1}$, \thanks{Partially supported by NSFC (No. 11771341). E-mail: chongzhi@whu.edu.cn}
\qquad  Duanzhi Zhang$^{2}$ \thanks{Partially supported by NSFC (No.11422103, 11171341, 11271200) and LPMC of Nankai
University. E-mail: zhangdz@nankai.edu.cn}\\ \\
$^{1}$ School of Mathematics and Statistics, Wuhan University,
\\Wuhan 430072, Hubei, P. R. China \\
$^{2}$ School of Mathematical Sciences and LPMC, Nankai University,\\ Tianjin 300071, P. R. China }
\date{}

\begin{document}

\maketitle

\begin{abstract}
{\it  Let $\Sigma$ be a compact convex
hypersurface in ${\bf R}^{2n}$ which is P-cyclic symmetric, i.e., $x\in \Sigma$ implies $Px\in\Sigma$ with P being a $2n\times2n$ symplectic orthogonal
matrix and $P^k=I_{2n}$, where $n, k\geq2$, $ker(P-I_{2n})=0$. In this paper, we first generalize Ekeland index theory for
periodic solutions of convex Hamiltonian system to a index theory with P boundary value condition and study its relationship with Maslov P-index theory,
then we use index theory to prove the existence of elliptic and non-hyperbolic closed characteristics on
compact convex P-cyclic symmetric hypersurfaces in ${\bf R}^{2n}$ for a broad class of symplectic orthogonal matrix P.}
\end{abstract}



\renewcommand{\theequation}{\thesection.\arabic{equation}}
\renewcommand{\thefigure}{\thesection.\arabic{figure}}

\setcounter{equation}{0}
\section{Introduction and main results}
Let $\Sigma$ be a $C^2$ compact hypersurface in ${\bf R}^{2n}$,
bounding a strictly convex compact set $U$ with non-empty interior, where $n\geq2$.
We denote the set of all such hypersurfaces by $\mathcal{H}(2n)$.
Without loss of generality, we suppose $U$ contains the origin. We
consider closed characteristics $(\tau,y)$ on $\Sigma$, which are
solutions of the following problem \be
\left\{\matrix{\dot{y}(t)=JN_\Sigma(y(t)),~y(t)\in \Sg,~\forall~t\in
{\bf R}, \cr
               y(\tau)=y(0), ~~~~~~~~~~~~~~~~~~~~~~~~~~~~~~~~~\cr }\right. \lb{1.1}\ee
where $J=\left(
              \begin{array}{cc}
                0 & -I_{n} \\
                I_{n} & 0 \\
              \end{array}
            \right)
$, $I_{n}$ is the identity matrix in ${\bf R}^n$ and
$\mathit{N}_\Sigma(y)$ is the outward normal unit vector  of
$\Sigma$ at $y$ normalized by the condition
$\mathit{N}_\Sigma(y)\cdot y=1$. Here $a\cdot b$ denotes the
standard inner product of $a, b \in {\bf R}^{2n}$. A closed
characteristic $(\tau, y)$ is {\it prime} if $\tau$ is the minimal
period of $y$. Two closed characteristics $(\tau,x)$ and
$(\sigma,y)$ are {\it geometrically distinct}, if $x({\bf
R})\not=y({\bf R})$. We denote by $\mathcal{J}(\Sigma)$ the set of all closed characteristics
$(\tau,y)$ on $\Sigma$ with $\tau$ being the minimal period of $y$.
Let $P$ be a $2n\times2n$ symplectic orthogonal
matrix and $P^k=I_{2n}$, where $k\geq2$. As in \cite{Zha1}, we call $\Sigma$ {\it P-cyclic symmetric}
if $P\Sigma=\Sigma$, i.e., $x\in \Sigma$ implies $Px\in\Sigma$. We denote by $\mathcal {H}_P(2n)$ the set of all P-cyclic symmetric hypersurfaces
in $\mathcal {H}(2n)$. A closed characteristic $(\tau, y)$ on $\Sigma \in
\mathcal {H}_P(2n)$ is {\it P-cyclic symmetric} if $y({\bf
R})=Py({\bf R})$, cf. Proposition 1 of \cite{Zha1}. In this paper, we further assume $ker(P-I_{2n})=0$.

Let $j: {\bf R}^{2n}\to {\bf R}$ be the gauge function of $\Sigma$,
i.e., $j(\lambda x)=\lambda$ for $x\in\Sigma$ and $\lambda\geq 0$,
then $j\in C^{2}({\bf R}^{2n}\setminus\{0\},{\bf R})\cap C^{1}({\bf
R}^{2n},{\bf R})$ and $\Sigma=j^{-1}(1)$. Fix a constant $\alpha
\in(1, 2)$ and define the Hamiltonian $H: {\bf
R}^{2n}\to[0,+\infty)$ by\[ H(x) :=j(x)^{\alpha}\] Then
$H\in C^{2}({\bf R}^{2n}\setminus\{0\},{\bf R})\cap
C^{0}({\bf R}^{2n},{\bf R})$ is convex and
$\Sigma=H^{-1}(1)$. It is well known that the problem
$(1.1)$ is equivalent to the following given energy problem of the
Hamiltonian system
\bea\left\{\begin{array}{ll}
\dot{y}(t)=JH^\prime(y(t)), H(y(t))=1,~\forall~t\in {\bf R},\\
y(\tau)=y(0).
\end{array}\right. \lb{1.2}\eea
Denote by $\mathcal{J}(\Sigma,\alpha)$ the set of all solutions
$(\tau,y)$ of the problem $(1.2)$, where $\tau$ is the minimal
period of $y$. Note that elements
in $\mathcal{J}(\Sigma)$ and $\mathcal{J}(\Sigma,\alpha)$ are in one
to one correspondence with each other. Let $(\tau,y)
\in\mathcal{J}(\Sigma,\alpha)$. We call the fundamental solution
$\gamma_y: [0,\tau]\to Sp(2n)$ with $\gamma_y(0)=I_{2n}$ of the
linearized Hamiltonian system\bea
\dot{z}(t)=JH^{\prime\prime}(y(t))z(t),~\forall~t \in {\bf
R}.\lb{1.3}\eea the {\it associated symplectic path} of $(\tau,y)$. The
eigenvalue of $\gamma_y(\tau)$ are called  {\it Floquet multipliers}
of $(\tau,y)$. By Proposition 1.6.13 of \cite{Eke2}, the Floquet
multipliers with their multiplicities and Krein type numbers of
$(\tau,y)\in\mathcal{J}(\Sigma,\alpha)$ do not depend on the
particular choice of the Hamiltonian function in $(1.2)$. As in
Chapter 15 of \cite{Lon4}, for any symplectic matrix M, we define the elliptic height $e(M)$ of
M by the total algebraic multiplicity of all eigenvalues of M on the unit circle
${\bf U}$ in the complex plane ${\bf C}$. And for any $(\tau,y)\in\mathcal{J}(\Sigma,\alpha)$
we define $e(\tau,y)=e(\gamma_y(\tau))$, and call $(\tau,y)$ {\it elliptic} or {\it hyperbolic} if
$e(\tau,y)=2n$ or $e(\tau,y)=2$, respectively.

The study on closed characteristics in the global sense started in 1978, when the existence of at least one
closed characteristic was first established on any compact star-shaped hypersurface by P. Rabinowitz in \cite{Rab1}
and on any compact convex hypersurface by A. Weinstein in \cite{Wei1} independently, since then the existence of
multiple closed characteristics on $\Sg\in\H(2n)$ has been deeply studied by many mathematicians, for
example, studies in \cite{EkL1}, \cite{EkH1}, \cite{Szu1}, \cite{HWZ1}, \cite{LLZ1}, \cite{LoZ1}, \cite{WHL1} and \cite{Wan2}.
There is a long standing conjecture on the stability of closed
characteristics on compact convex hypersurfaces in ${\bf
R}^{2n}$: for every $\Sigma\in\mathcal {H}(2n)$, there exists an elliptic $(\tau,y)\in\mathcal{J}(\Sigma)$.
cf., Page 235 of \cite{Eke2}. Ekeland proved in
\cite{Eke1} of 1986 the existence of at least one elliptic closed
characteristic on $\Sg$ provided $\Sg\in\H(2n)$ is
$\sqrt{2}$-pinched. In \cite{DDE1} of 1992, Dell'Antonio, D'Onofrio
and Ekeland proved the existence of at least one elliptic closed
characteristic on $\Sg$ provided $\Sg\in\mathcal{H}_P(2n)$ if $P=-I_{2n}$. In \cite{LoZ1} of 2002, Long and
Zhu proved when $^\#{\J}(\Sg)<+\infty$, there exists at
least one elliptic closed characteristic. For more results on the stability of closed
characteristics on compact convex hypersurfaces, we refer readers to \cite{HuO, Lon3, Lon4, LoZ1, Wan1, WHL1} and the reference therein.

In this paper, we prove the above conjecture for a broad class of compact convex P-cyclic symmetric hypersurfaces.
For the studies about closed characteristics on compact convex P-cyclic symmetric hypersurfaces, one can also refer to \cite{DoL1, DoL2, Liu1, LiZ1, LiZ2, Zha1}.

{\bf Theorem 1.1.} {\it Assume $\Sigma\in\mathcal {H}_P(2n)$ and $P$ has the form $P=R(\theta_1)\diamond\cdots \diamond R(\theta_n)$ for $\theta_i\in(0, \pi]$
, where $R(\th)=\left(\matrix{\cos\th & -\sin\th\cr
        \sin\th  & \cos\th\cr}\right)$. Then there exist at least one elliptic closed characteristic on $\Sigma$.}

When we weaken the condition on P, the existence of one non-hyperbolic closed characteristic can be obtained:

{\bf Theorem 1.2.} {\it Assume $\Sigma\in\mathcal {H}_P(2n)$ and $P$ has the form $P=R(\theta_1)\diamond\cdots \diamond R(\theta_n)$ for $\theta_i\in(0, 2\pi)$
with $^\#\{i\mid \theta_i\in(0, \pi]\}-^\#\{i\mid \theta_i\in(\pi, 2\pi)\}\geq 2$, where $R(\th)=\left(\matrix{\cos\th & -\sin\th\cr
        \sin\th  & \cos\th\cr}\right)$. Then there exist at least one non-hyperbolic closed characteristic on $\Sigma$.}

This paper is arranged as follows. In Section 2, we first generalize Ekeland index theory for
periodic solutions of convex Hamiltonian system to a index theory with P boundary value condition. In Section 3,
we recall briefly the Maslov P-index theory for symplectic paths and study its relationship with
Ekeland P-index theory. In Section 4, we establish a variational structure for closed characteristics and use
index theory to prove Theorem 1.1 and Theorem 1.2.

\setcounter{equation}{0}
\section{Ekeland P-index theory for positive definite Hamiltonian system}
In this section we will offer a slight generalization of the contents of Section 1.4 in \cite{Eke2}.

Let $A(t)$ be a symmetric and positive definite $2n\times 2n$ real matrix with depending continuously on $t\in[0,+\infty)$.
Then we consider the following quadratic form:
\begin{equation}\label{2.1}
q_{s}(u,u)=\frac{1}{2}\int_0^s [(Ju,\Pi_{s}u)+(B(t)Ju,Ju)]\,dt, \qquad\forall u\in L^{2}(0,s),
\end{equation}
where $B(t)=A(t)^{-1}$, $L^{2}(0,s)=L^{2}((0,s),\mathbf{R}^{2n})$ and $\Pi_{s}:L^{2}(0,s)\rightarrow L^{2}(0,s)$ is defined by
\begin{equation}\label{2.2}
(\Pi_{s}u)(t)=x(t),
\end{equation}
\begin{equation}\label{2.3}
x(t)=\int_0^t u(\tau)\,d\tau+(P-I)^{-1}\int_0^s u(\tau)\,d\tau.
\end{equation}
Here $P$ is an orthogonal symplectic matrix satisfying $P^k=I_{2n}$ for  integer $k\in[2,+\infty)$ and
$\ker(P-I_{2n})=0$. Note that $x(s)=Px(0)$. In the following, we denote $L^2=L^2(0,s)$ for simplicity.

{\bf Lemma 2.1.} {\it $\Pi_{s}$ is a compact operator  from $L^{2}(0,s)$ into itself. Moreover $\Pi_{s}$ is antisymmetric. }

{\bf Proof.} $\Pi_{s}$ sends $L^{2}(0,s)$ into $W^{1,2}(0,s)$ and the identity map from $W^{1,2}(0,s)$ to $L^{2}(0,s)$ is compact
by the Rellich-Kondrachov theorem.

To check that it is antisymmetric, we just integrate by parts:
\[
\int^s_0 (\Pi_{s}u,v)\,dt=-\int^s_0 (u,\Pi_{s}v)\,dt+(\Pi_{s}u,\Pi_{s}v)\,|{^s_0},
\]
and the last term vanishes since:
\begin{eqnarray*}
(\Pi_{s}u,\Pi_{s}v)\,|{^s_0}&=&((\Pi_{s}u)(s),(\Pi_{s}v)(s))-((\Pi_{s}u)(0),(\Pi_{s}v)(0))  \\
&=&(x(s),y(s))-(x(0),y(0))  \\
&=&(Px(0),Py(0))-(x(0),y(0))\\
&=&(P^TPx(0),y(0))-(x(0),y(0)) \\
&=&0,
\end{eqnarray*}
where $y(t)=(\Pi_sv)(t)$.\hfill\hb

Thus, by Lemma 2.1, we have $J\Pi_{s}:L^{2}(0,s)\rightarrow L^{2}(0,s)$ is self-adjoint and compact.

{\bf Lemma 2.2.} {\it For any symmetric and positive definite $2n \times 2n$ real matrix $A(t)$ continuous in $t\in [0,+\infty)$, there is a splitting:
\[
L^{2}(0,s)=E^{+}(A)\oplus E^{0}(A)\oplus E^{-}(A)
\]
such that:

(a)$E^{+}(A)$, $E^{0}(A)$, $E^{-}(A)$ are $q_{s}$-orthogonal,

(b)$q_{s}(u,u)>0$\quad $\forall u\in E^{+}(A)\backslash\{0\}$,

(c)$q_{s}(u,u)=0$\quad $\forall u\in E^{0}(A)$,

(d)$q_{s}(u,u)<0$\quad $\forall u\in E^{-}(A)\backslash\{0\}$,

(e)$E^{0}(A)$ and $E^{-}(A)$ are finite-dimensional.}

{\bf Proof.} Define a self-adjoint operator $\bar{B}:L^{2}\rightarrow L^{2}$ by
\[
(\bar{B}u,v)=\int_0^s (B(t)Ju(t), Jv (t))\,dt,\qquad \forall u,v\in L^2\equiv L^{2}(0,s).
\]
Since $A(t)$ is symmetric, positive definite, and continuous, we can find positive constants $a$ and $b$ such that
\[
a(x,x)\geq(B(t)x,x)\geq b(x,x),\quad \forall x\in\mathbf{R}^{2n},t\in[0,s].
\]
This yields
\[
a\|u\|^2\geq(\bar{B}u,u)\geq b\|u\|^2,\quad \forall u\in L^2(0,s).
\]
Hence, the Lax-Milgram theorem tells us that $\bar{B}$ is an isomorphism, and $(\bar{B}u,v)$ defines a Hilbert space structure on $L^{2}(0,s)$.
Endowing $L^{2}(0,s)$ with  the interior product $(\bar{B}u,v)$, and noticing that $\bar{B}^{-1}J\Pi_{s}$ is self-adjoint,
and applying to $\bar{B}^{-1}J\Pi_{s}$ the spectral theory of compact self-adjoint operators on a Hilbert space,
we know that there is a basis $\{e_j\}_{j\in \mathds{N}}$ of $L^{2}(0,s)$, and a sequence $\lambda_j\rightarrow 0$ in $\mathbf{R}$
as $j\rightarrow +\infty$ such that
\begin{eqnarray*}
(\bar{B}e_{i}, e_{j})&=&\delta_{ij}, \\
\bar{B}^{-1}J\Pi_{s}e_{j}&=&\lambda_{j}e_{j}.
\end{eqnarray*}
Therefore, for any $u=\sum\nolimits_{j=1}^{\infty}c_{j}e_{j}\in L^{2}(0,s)$, by the definition of formula (\ref{2.1}), we obain
\begin{eqnarray*}
q_{s}(u,u)&=&-\frac{1}{2}(J\Pi_{s}u,u)+\frac{1}{2}(\bar{B}u, u) \\
&=& -\frac{1}{2}\sum_{j=1}^\infty\lambda_jc^2_j+\frac{1}{2}\sum_{j=1}^\infty c^2_j \\
&=& \frac{1}{2}\sum_{j=1}^\infty(1-\lambda_j)c^2_j
\end{eqnarray*}

Since $\lambda_j\rightarrow 0$  as $j\rightarrow +\infty$, all the coefficients $(1-\lambda_j)$ are positive except  a finite number. Thus the result, with:
\begin{eqnarray*}
E^{+}(A)&=&\left\{\sum c_{j}e_{j}\;|\;c_{j}=0 \quad{\rm if}\quad 1-\lambda_{j}\leq0 \right\},\\
E^{0}(A)&=&\left\{\sum c_{j}e_{j}\;|\;c_{j}=0 \quad{\rm if}\quad 1-\lambda_{j}\neq0 \right\},\\
E^{-}(A)&=&\left\{\sum c_{j}e_{j}\;|\;c_{j}=0 \quad{\rm if}\quad 1-\lambda_{j}\geq0 \right\}.
\end{eqnarray*}

{\bf Definition 2.3.} {\it For any symmetric and positive definite $2n \times 2n$ real matrix $A(t)$ continuous in $t\in [0,s]$, we define}
\bea
\nu^{E}_{P}(A)=\dim E^{0}(A),\qquad i^{E}_{P}(A)=\dim E^{-}(A).
\nn\eea

{\bf Proposition 2.4.} {\it For any symmetric and positive definite $2n \times 2n$ real matrix $A(t)$ continuous in $t\in [0,s]$, $\nu^{E}_{P}(A)$ is
the number of linearly independent solutions of the following problem:
\[
\left\{
\begin{array}{l}
\dot{x}(t)=JA(t)x\\
 x(s)=Px(0).
\end{array}
\right.
\]
In other words, $\nu^{E}_{P}(A)=\dim\ker(\gamma_A(s)-P)$, where $\gamma=\gamma_A(t)$ is the fundamental solution of $\dot{x}(t)=JA(t)x$
with $\gamma(0)=I_{2n}$.
}

{\bf Proof.} For any $u\in E^{0}(A)$, according to the definition of formula (\ref {2.1}) and Lemma 2.2, we obtain
\begin{eqnarray*}
q_{s}(u,v)&=&\frac{1}{2}\int_0^s [(Ju,\Pi_{s}v)+(B(t)Ju,Jv)]\,dt \\
&=&\frac{1}{2}\int_0^s (-J\Pi_{s}u-JB(t)Ju,v)\,dt=0, \qquad\forall v\in L^{2}(0,s).
\end{eqnarray*}
The kernel of $q_{s}$ consists of all $u\in L^{2}(0,s)$ such that this interior product vanishes for all $v \in L^{2}(0,s)$. Thus, we have
\[
-J\Pi_{s}u-JB(t)Ju=0,
\]
which yields
\begin{equation}\label{2.4}
JA(t)\Pi_{s}u=u.
\end{equation}
Now define $x=\Pi_{s}u$. We obtain $u=\dot{x}$, $x(s)=Px(0)$ and formula (\ref{2.4}) reads as
\[
\dot{x}=JA(t)x \quad {\rm for} \quad t \in (0,s).
\]
Moreover, we get $x(t)=\gamma_{A}(t)c$, where $c\in \mathbf{R}^{2n}$ satisfies
\[
\gamma_{A}(s)c=x(s)=Px(0)=P\gamma_{A}(0)c=PI_{2n}c=Pc.
\]
This yields,
\begin{equation}\label{2.5}
(\gamma_{A}(s)-P)c=0.
\end{equation}
Hence we obtain
\[
E^{0}(A)\cong\left\{   c\in \mathbf{R}^{2n}|(\gamma_{A}(s)-P)c=0   \right\}=\ker(\gamma_A(s)-P).
\]
So
\[
\nu_P^E(A)=\dim E^0(A)=\dim\ker(\gamma_A(s)-P).
\]

{\bf Proposition 2.5.} {\it For any symmetric and positive definite $2n \times 2n$ real matrix $A(t)$ continuous in $t\in [0,s]$, we have
\begin{equation}\label{2.6}
i^E_P(A)=\sum_{0<\sigma<s}\nu^E_P(A_\sigma),
\end{equation}
where $A_\sigma=A|_{[0,\sigma]}$.}

{\bf Proof.} The proof proceeds through five steps.

{\bf Step 1.} When $\sigma>0$ sufficiently small, we have $i^E_P(A_\sigma)=0$.

In fact, by formulas (\ref {2.2})-(\ref {2.3}) and Cauchy-Schwarz inequality, we have
\begin{eqnarray*}
|x(t)|&\leq&\int_0^t |u(\tau)|\,d\tau+\|(P-I)^{-1}\|\int_0^\sigma |u(\tau)|\,d\tau \\
&\leq&C\int_0^\sigma |u(\tau)|\,d\tau\\
&\leq&C\left(\int_0^\sigma |u(\tau)|^2\,d\tau\right)^{\frac{1}{2}}\left(\int_0^\sigma 1^2\,d\tau\right)^{\frac{1}{2}}\\
&=&C\sigma^{\frac{1}{2}}\|u\|,
\end{eqnarray*}
where $C=1+\|(P-I)^{-1}\|>1$, and $\|\cdot\|$ is the $L^2$-norm. Hence,
\begin{eqnarray*}
\|x\|&=&\left(\int_0^\sigma |x(t)|^2\,dt\right)^{\frac{1}{2}} \\
&\leq&\left(\int_0^\sigma \left(C\sigma^{\frac{1}{2}}\|u\|\right)^2\,dt\right)^{\frac{1}{2}}\\
&=&\left(C^2\sigma^2\|u\|^2\right)^\frac{1}{2}\\
&=&C\sigma\|u\|.
\end{eqnarray*}
This yields
\[
\|\Pi_{\sigma}u\|=\|x\|\leq C\sigma\|u\|,\quad \forall u\in L^2(0,\sigma).
\]
In addition, since $A(t)$ is symmetric, positive definite, and continuous, we can find a positive constant $b$ such that
\[
(B(t)x,x)\geq b(x,x),\quad \forall x\in\mathbf{R}^{2n},t\in[0,\sigma].
\]
Applying Cauchy-Schwarz to formula (\ref{2.1}), we have
\begin{eqnarray*}
q_{\sigma}(u,u)&=&\frac{1}{2}\int_0^\sigma [(Ju,\Pi_{\sigma}u)+(B(t)Ju,Ju)]\,dt\\
&\geq&\frac{1}{2}(-\|u\|\,\|\Pi_{\sigma}u\|+b\|u\|^2)\\
&\geq&\frac{1}{2}(b-C\sigma)\|u\|^2.
\end{eqnarray*}
So $q_{\sigma}$ is positive definite for $\sigma<\frac{b}{C}$.

{\bf Step 2.} We claim that there are only finitely many points $\sigma$ with $\nu^E_P(A_\sigma)\neq 0$ in any bounded interval $[0,s]$.

Now argue by contradiction. In fact, if not, by (\ref{2.5}) there exist $\lambda_j\in[0,s]$ and $\xi_j\in\mathbf{R}^{2n}\backslash\{0\}$ with $|\xi_j|=1$
such that
\begin{equation}\label{2.7}
\gamma_A(\lambda_{j})\xi_j=P\xi_j,\qquad {\rm for}\; j=1,2,\ldots
\end{equation}
Without loss of generality, we may assume that $\lambda_{j}\rightarrow \lambda$ and $\xi_j\rightarrow \xi$ as $j\rightarrow +\infty$. This yields
\begin{equation}\label{2.8}
\gamma_A(\lambda)\xi=P\xi,
\end{equation}
\begin{equation}\label{2.9}
(\gamma_A(\lambda_{j})-P)(\xi_j-\xi)=(\gamma_A(\lambda)-\gamma_A(\lambda_j))\xi.
\end{equation}
As we know $\gamma_A(\lambda_j)$ and $P$ are symplectic, then we have $\gamma_A(\lambda_j)^T J=J\gamma_A(\lambda_j)^{-1}$ and $P^TJP=J$.
Lastly, by formula (\ref {2.7}) we get $\gamma_A(\lambda_j)^{-1} P\xi_j=\xi_j$. Hence,
\begin{eqnarray*}
(\gamma_A(\lambda_j)(\xi_j-\xi),JP\xi_j)&=&(\xi_j-\xi,\gamma_A(\lambda_j)^T JP\xi_j) \\
  &=&(\xi_j-\xi,J\gamma_A(\lambda_j)^{-1} P\xi_j)  \\
  &=&(\xi_j-\xi,J\xi_j) \\
  &=&(\xi_j-\xi,P^TJP\xi_j)\\
  &=&(P(\xi_j-\xi),JP\xi_j).
\end{eqnarray*}
Thus $((\gamma_A(\lambda_j)-P)(\xi_j-\xi),JP\xi_j)=0$. What's more, by formula (\ref {2.9}) we have
\[
((\gamma_A(\lambda)-\gamma_A(\lambda_j))\xi,JP\xi_j)=0.
\]
This yields
\[
0=\lim_{j\rightarrow+\infty}(\frac{\gamma_A(\lambda)-\gamma_A(\lambda_j)}{\lambda-\lambda_j}\xi,JP\xi_j)=(\dot{\gamma}_A(\lambda)\xi,JP\xi).
\]
In addition, as we know $\dot{\gamma}_A(\lambda)=JA(\lambda)\gamma_A(\lambda)$ and from formula (\ref{2.8}), we obtain
\begin{eqnarray*}
0&=&(\dot{\gamma}_A(\lambda)\xi,JP\xi) \\
&=&(JA(\lambda)\gamma_A(\lambda)\xi,JP\xi)\\
&=&(JA(\lambda)P\xi,JP\xi)\\
&=&(A(\lambda)P\xi,P\xi),
\end{eqnarray*}
which contradicts to the fact that $A(\lambda)$ is positive definite.

{\bf Step 3.} If $\sigma_1<\sigma_2$, there hold
\begin{equation}\label{2.10}
i_P^E(A_{\sigma_1})\leq i_P^E(A_{\sigma_2}),
\end{equation}
\begin{equation}\label{2.11}
i_P^E(A_{\sigma_1})+\nu_P^E(A_{\sigma_1})\leq i_P^E(A_{\sigma_2}).
\end{equation}
In fact, we define a map $\theta:L^2(0,\sigma_1)\rightarrow L^2(0,\sigma_2)$ by
\[ (\theta u)(t) = \left\{\matrix{
            u(t), & \quad {\rm if}\;0\leq t\leq \sigma_1, \cr
            0, & \quad {\rm if}\; \sigma_1< t\leq \sigma_2. \cr}\right. \]
Clearly, for any $u \in L^2(0,\sigma_1)$ we get
\[
q_{\sigma_2}(\theta u,\theta u)=q_{\sigma_1}( u, u).
\]
Therefore,
\[
q_{\sigma_2}( u, u)<0, \quad \forall u\in\theta(E^-(A_{\sigma_1}))\setminus\{0\}.
\]
This yields
\[
i^E_P(A_{\sigma_2})\geq\dim(\theta(E^-(A_{\sigma¡ª¡ª1})))=i^E_P(A_{\sigma_1}).
\]
Hence the proof of (\ref{2.10}) is done. Similarly, we obtain
\begin{equation}\label{2.12}
i_P^E(A_{\sigma_1})+\nu_P^E(A_{\sigma_1})\leq i_P^E(A_{\sigma_2})+\nu_P^E(A_{\sigma_2}).
\end{equation}
Moreover, formula (\ref{2.11}) follows from (\ref{2.10}) when $\nu_P^E(A_{\sigma_1})=0$. On the other hand, if $\nu_P^E(A_{\sigma_1})\neq0$,
then from Step 2 we have $\nu_P^E(A_{\sigma_1^+})=0$. Hence by (\ref{2.12}) and let $\sigma_2=\sigma_1^+$, we have
\[
i_P^E(A_{\sigma_1})+\nu_P^E(A_{\sigma_1})\leq i_P^E(A_{\sigma_1^+}).
\]
However, if $\sigma_2>\sigma_1$ then $\sigma_2\geq\sigma_1^+$, by (\ref{2.10}), we get
\[
i_P^E(A_{\sigma_1^+})\leq i_P^E(A_{\sigma_2}).
\]
Therefore,
\[
i_P^E(A_{\sigma_1})+\nu_P^E(A_{\sigma_1})\leq i_P^E(A_{\sigma_1^+})\leq i_P^E(A_{\sigma_2}).
\]

{\bf Step 4.} The function $s\rightarrow i^E_P(A_s)$ is left continuous, i.e., $i^E_P(A_s)=i^E_P(A_{s^-})$

In fact, let $u_1(t)=u(st)$, by formulas (\ref{2.2})-(\ref{2.3}), we have,
\[
(\Pi_{1}u_1)(t)=\int_0^t u(s\tau)\,d\tau+(P-I)^{-1}\int_0^1 u(s\tau)\,d\tau.
\]
Let $\alpha=st$. Then we calculate $(\Pi_{s}u)(\alpha)$ as follows
\bea (\Pi_{s}u)(\alpha)&=&\int_0^\alpha u(\tau)\,d\tau+(P-I)^{-1}\int_0^s u(\tau)\,d\tau\nn\\
&=&s\int_0^t u(s\tau)\,d\tau+s(P-I)^{-1}\int_0^1 u(s\tau)\,d\tau\nn \\
&=&s(\Pi_{1}u_1)(t).\label{2.13}\eea
Define a map $p:L^2(0,s)\rightarrow L^2(0,1)$ by $(pu)(t)=u(st)=u_1(t)$. And define a quadratic form on $L^2(0,1)$ by:
\begin{equation}\label{2.14}
q_{s}^1(u,u)=\frac{1}{2}\int_0^1 [s(Ju,\Pi_{1}u)+(B(st)Ju,Ju)]\,dt.
\end{equation}
Next, we will prove that $q_{s}(u,u)=sq_{s}^1(pu,pu)$. In fact, by formulas (\ref{2.1}), (\ref{2.13}) and (\ref{2.14}), we obtain
\begin{eqnarray*}
sq_{s}^1(pu,pu)&=& \frac{s}{2}\int_0^1 [s(Jpu,\Pi_{1}pu)+(B(st)Jpu,Jpu)]\,dt \\
&=& \frac{1}{2}\int_0^1 [(Ju(st),s(\Pi_{1}u_1)(t))+(B(st)Ju(st),Ju(st))]\,dst\\
&=& \frac{1}{2}\int_0^s [(Ju(\alpha),(\Pi_{s}u)(\alpha))+(B(\alpha)Ju(\alpha),Ju(\alpha))]\,d\alpha\\
&=& q_{s}(u,u).
\end{eqnarray*}
Hence, for any fixed $s_0$, let $E_1=p(E^-(A_{s_0}))$, then we have
\[
i^E_P(A_{s_0})=\dim E_1,
\]
\[
q^1_{s_0}(u,u)<0,\quad \forall u\in E_1\backslash\{0\}.
\]
Since $q^1_{s}$ depends continuously on $s$ in formula (\ref{2.14}), as $s\rightarrow s_0$ we obtain
\[
q^1_{s}(u,u)<0,\quad \forall u\in E_1\backslash\{0\}.
\]
This yields
\[
i^E_P(A_{s_0})\leq i^E_P(A_{s})\Rightarrow i^E_P(A_{s})\leq i^E_P(A_{s^-}).
\]
The converse inequality holds by (\ref{2.10}) of Step 3. Therefore $i^E_P(A_s)=i^E_P(A_{s^-})$.

{\bf Step 5.} For any $\sigma\in[0,s)$, there holds
\begin{equation}\label{2.15}
i_P^E(A_{\sigma^+})=i^E_P(A_\sigma)+\nu^E_P(A_\sigma).
\end{equation}
Moreover, $i^E_P(A_\sigma)$ is continuous at the point $\sigma\in(0,s)$ with $\nu^E_P(A_\sigma)=0$.

In fact, by Step 4 we know $i^E_P(A_s)$ and $\nu^E_P(A_s)$ are also the index and nullity of $q^1_{s}$
which is defined in formula (\ref{2.14}) on $L^2(0,1)$. For the sake of convenience, we just consider $q^1_{s}$. Denote by
\[
(B^1(s)u,v)=\int_0^1 (B(st)Ju, Jv)\,dt,\qquad \forall u,v\in L^{2}(0,1).
\]
Arguing as in Lemma 2.2, we know that there is a basis $\{e_j^s\}_{j\in \mathds{N}}$ of $L^{2}(0,1)$,
and a sequence $\lambda_j^s\rightarrow 0$ in $\mathbf{R}$ as $j\rightarrow +\infty$ such that
\begin{eqnarray*}
(B^1(s)e_i^s,e_j^s)&=&\delta_{ij},    \\
(J\Pi_{1}e_j^s,u)&=&\lambda_j^s(B^1(s)e^s_j,u),\qquad \forall u\in L^2(0,1).
\end{eqnarray*}
Therefore, for any $u=\sum\nolimits_{j=1}^{\infty}\xi_je_j^s\in L^{2}(0,1)$, by the definition of formula (\ref{2.14}), we obtain
\[
q^1_{s}(u,u)=\frac{1}{2}\sum_{j=1}^\infty(1-s\lambda_j^s)\xi^2_j.
\]
For any fixed $\sigma>0$, we set $i^E_P(A_{\sigma^+})=K$. This means that there is a $\sigma'>\sigma$ such that $i^E_P(A_s)=K$ for
all $s\in(\sigma,\sigma')$. Thus for any $s\in(\sigma,\sigma')$ we obtain
\[
1-s\lambda^s_j<0\quad {\rm for}\; 1\leq j\leq K.
\]
Fix $j\leq K$. As we know the $e^s_j$ and $\lambda^s_j=(J\Pi_{1} e^s_j,e^s_j)$ are bounded with $\lambda^s_j>\frac{1}{s}>\frac{1}{\sigma'}$.
Then there exist $\{e_j^{s(l)}\}$ and $\{\lambda_j^{s(l)}\}$ such that $e_j^{s(l)}\rightarrow e_j$ in $L^2(0,1)$
and $\lambda_j^{s(l)}\rightarrow \lambda_j$, $s(l)\rightarrow \sigma$ in $\mathbf{R}$ as $l\rightarrow +\infty$.
Thus we obtain
\begin{eqnarray*}
(B^1(\sigma)e_i,e_j) &=&\delta_{ij}, \quad {\rm for}\;i,j=1,2,\ldots,K, \\
(J\Pi_{1}e_j,u)&=&\lambda_j(B^1(\sigma) e_j,u), \quad \forall u\in L^2(0,1),j=1,2,\ldots,K,  \\
   1-\sigma\lambda_j&\leq&0, \quad {\rm for}\;j=1,2,\ldots,K.
\end{eqnarray*}
Therefore, for any $u=\sum\nolimits_{j=1}^{K}\xi_je_j\in L^{2}(0,1)$, by the definition of formula (\ref{2.14}), we obtain
\[
q^1_{\sigma}(u,u)=\frac{1}{2}\sum_{j=1}^K(1-\sigma\lambda_j)\xi^2_j\leq0.
\]
So that,
\[
i^E_P(A_{\sigma^+})=K\leq i^E_P(A_\sigma)+\nu^E_P(A_\sigma).
\]
The converse inequality holds by (\ref{2.11}) of  Step 3. Therefore $i_P^E(A_{\sigma^+})=i^E_P(A_\sigma)+\nu^E_P(A_\sigma)$.
Moreover, if $\nu^E_P(A_\sigma)=0$, by Step 4 and formula (\ref{2.15}), we obtain
\[
i^E_P(A_{\sigma^-})=i^E_P(A_\sigma)=i_P^E(A_{\sigma^+}).
\]
Hence, we get the results.

{\bf Proof of Proposition 2.5.} The function $\sigma\rightarrow i^E_P(A_\sigma)$ is integer-valued, left continuous and non-decreasing on $(0,+\infty)$.
Its value at any point $s$ must be equal to the sum of the jumps it incurred in $(0,s)$. By Step 5,
this is precisely the sum of the $\nu^E_P(A_\sigma)$ with $0<\sigma<s$.\hfill\hb

\setcounter{equation}{0}
\section{Relationship between Ekeland P-index theory with Maslov P-index theory }
In this section, we recall briefly the Maslov P-index theory for symplectic paths and study its relationship with
Ekeland P-index theory. Note that the Maslov P-index theory
for a symplectic path was first studied by Y. Dong and C. Liu in \cite{Dong, LiuC} independently for any symplectic matrix P with different
treatment. The Maslov P-index theory was generalized in \cite{LT1} to the Maslov
$(P, \omega)$-index theory for any $P\in Sp(2n)$ and all $\omega\in \U$. The iteration theory of
$(P, \omega)$-index theory was studied in \cite{LT2}. When $\omega= 1$, the Maslov $(P, \omega)$-index theory
coincides with the Maslov P-index theory.

As usual, the symplectic group $\Sp(2n)$ is defined by
$$ \Sp(2n) = \{M\in {\rm GL}(2n,\R)\,|\,M^TJM=J\}, $$
whose topology is induced from that of $\R^{4n^2}$. For $\tau>0$ we are interested
in paths in $\Sp(2n)$:
$$ \P_{\tau}(2n) = \{\ga\in C([0,\tau],\Sp(2n))\,|\,\ga(0)=I_{2n}\}.$$
We consider this path-space equipped
with the $C^0$-topology.
For any $\om\in\U$ the following codimension $1$ hypersurface in $\Sp(2n)$ is
defined in \cite{Lon2}:
$$ \Sp(2n)_{\om}^0 = \{M\in\Sp(2n)\,|\, \det(M-\om I_{2n}))=0\}.  $$
For any $M\in \Sp(2n)_{\om}^0$, we define a co-orientation of $\Sp(2n)_{\om}^0$
at $M$ by the positive direction $\frac{d}{dt}Me^{t J}|_{t=0}$.  Let
\bea
\Sp(2n)_{\om}^{\ast} &=& \Sp(2n)\bs \Sp(2n)_{\om}^0,   \nn\\
\P_{\tau,\om}^{\ast}(2n) &=&
      \{\ga\in\P_{\tau}(2n)\,|\,\ga(\tau)\in\Sp(2n)_{\om}^{\ast}\}, \nn\\
\P_{\tau,\om}^0(2n) &=& \P_{\tau}(2n)\bs  \P_{\tau,\om}^{\ast}(2n).  \nn\eea
For any two continuous arcs $\xi$ and $\eta:[0,\tau]\to\Sp(2n)$ with
$\xi(\tau)=\eta(0)$, their concatenation is defined
as usual by
$$ \eta\ast\xi(t) = \left\{\matrix{
            \xi(2t), & \quad {\rm if}\;0\le t\le \tau/2, \cr
            \eta(2t-\tau), & \quad {\rm if}\; \tau/2\le t\le \tau. \cr}\right. $$
Given any two $2m_k\times 2m_k$ matrices of square block form
$M_k=\left(\matrix{A_k&B_k\cr
                                C_k&D_k\cr}\right)$ with $k=1, 2$,
as in \cite{Lon4}, the $\;\dm$-product of $M_1$ and $M_2$ is defined by
the following $2(m_1+m_2)\times 2(m_1+m_2)$ matrix $M_1\dm M_2$:
$$ M_1\dm M_2=\left(\matrix{A_1&  0&B_1&  0\cr
                               0&A_2&  0&B_2\cr
                             C_1&  0&D_1&  0\cr
                               0&C_2&  0&D_2\cr}\right). \nn$$  
Denote by $M^{\dm k}$ the $k$-fold $\dm$-product $M\dm\cdots\dm M$. Note
that the $\dm$-product of any two symplectic matrices is symplectic. For any two
paths $\ga_j\in\P_{\tau}(2n_j)$ with $j=0$ and $1$, let
$\ga_0\dm\ga_1(t)= \ga_0(t)\dm\ga_1(t)$ for all $t\in [0,\tau]$.

A special path $\xi_n$ is defined by
\be \xi_n(t) = \left(\matrix{2-\frac{t}{\tau} & 0 \cr
                                             0 &  (2-\frac{t}{\tau})^{-1}\cr}\right)^{\dm n}
        \qquad {\rm for}\;0\le t\le \tau.  \lb{3.1}\ee

{\bf Definition 3.1.} (cf. \cite{Lon2}, \cite{Lon4}) {\it For any $\om\in\U$ and
$M\in \Sp(2n)$, define
\be  \nu_{\om}(M)=\dim_{\C}\ker_{\C}(M - \om I_{2n}).  \lb{3.2}\ee
For any $\tau>0$ and $\ga\in \P_{\tau}(2n)$, define
\be  \nu_{\om}(\ga)= \nu_{\om}(\ga(\tau)).  \lb{3.3}\ee

If $\ga\in\P_{\tau,\om}^{\ast}(2n)$, define
\be i_{\om}(\ga) = [\Sp(2n)_{\om}^0: \ga\ast\xi_n],  \lb{3.4}\ee
where the right hand side of (\ref{3.4}) is the usual homotopy intersection
number, and the orientation of $\ga\ast\xi_n$ is its positive time direction under
homotopy with fixed end points.

If $\ga\in\P_{\tau,\om}^0(2n)$, we let $\mathcal{F}(\ga)$
be the set of all open neighborhoods of $\ga$ in $\P_{\tau}(2n)$, and define
\be i_{\om}(\ga) = \sup_{U\in\mathcal{F}(\ga)}\inf\{i_{\om}(\beta)\,|\,
                       \beta\in U\cap\P_{\tau,\om}^{\ast}(2n)\}.
               \lb{3.5}\ee
Then
$$ (i_{\om}(\ga), \nu_{\om}(\ga)) \in \Z\times \{0,1,\ldots,2n\}, $$
is called the index function of $\ga$ at $\om$. }

Note that when $\om=1$, this index theory was introduced by
C. Conley-E. Zehnder in \cite{CoZ1} for the non-degenerate case with $n\ge 2$,
Y. Long-E. Zehnder in \cite{LZe1} for the non-degenerate case with $n=1$,
and Y. Long in \cite{Lon1} and C. Viterbo in \cite{Vit1} independently for
the degenerate case. The case for general $\om\in\U$ was defined by Y. Long
in \cite{Lon2} in order to study the index iteration theory (cf. \cite{Lon4}
for more details and references).

For any $M\in\Sp(2n)$ and $\om\in\U$, the {\it splitting numbers} $S_M^{\pm}(\om)$
of $M$ at $\om$ are defined by
\be S_M^{\pm}(\om)
     = \lim_{\ep\to 0^+}i_{\om\exp(\pm\sqrt{-1}\ep)}(\ga) - i_{\om}(\ga),
   \lb{3.6}\ee
for any path $\ga\in\P_{\tau}(2n)$ satisfying $\ga(\tau)=M$, which is well defined by Lemma 9.1.5 of \cite{Lon4}.

{\bf Definition 3.2.} (cf. \cite{LT1}) {\it For any $P\in Sp(2n)$, $\omega\in \U$ and $\gamma \in \mathcal {P}_{\tau} (2n)$, the Maslov
$(P, \omega)$-index is defined by\bea i_\omega^P(\gamma)=i_\omega(P^{-1}\gamma*\xi)-i_\omega(\xi),\lb{3.7}\eea
where $\xi\in \mathcal {P}_{\tau} (2n)$ such that $\xi(\tau)=P^{-1}\gamma(0)=P^{-1}$, and $(P, \omega)$-nullity $\nu^P_
{\omega}(\gamma)$ is defined by} \bea \nu_\omega^P(\gamma)=dim_{{\bf C}} ker_{{\bf
C}}(\gamma(\tau) - \omega P),\lb{3.8}\eea

For any $M \in Sp(2n)$ and $\omega \in {\bf U}$, {\it the splitting
numbers} $_PS^{\pm}_M (\omega)$ of $M$ at $(P,\omega)$ are defined in Definition 2.4 of \cite{LT2} as follows
\begin{eqnarray}
_PS^{\pm}_M (\omega)=\lim_{\epsilon\rightarrow 0^+}{i^P_{\omega
\exp{(\pm\sqrt{-1}\epsilon)}}(\gamma)-i^P_{\omega}(\gamma)},
\lb{3.9}\end{eqnarray}
for any path $\gamma \in \mathcal {P}_{\tau} (2n)$ satisfying
$\gamma(\tau) = M$.

Let $\Omega^0(M)$ be the path connected component containing $M =
\gamma(\tau)$ of the set\begin{eqnarray} \Omega(M) = \{N \in Sp(2n)
\mid  \sigma(N)
\cap {\bf U} = \sigma(M) \cap {\bf U}~and~~~~\nn\\
\nu_{\lambda}(N)=\nu_{\lambda}(M),\forall \lambda \in \sigma(M) \cap
{\bf U}\}\lb{3.10}\end{eqnarray} Here $\Omega^0(M)$ is called the {\it
homotopy component} of $M$ in $Sp(2n)$.

In \cite{Lon2}-\cite{Lon4}, the following symplectic matrices were introduced as
basic normal forms: \begin{eqnarray}D(\lambda)=\left(
                                                    \begin{array}{cc}
                                                      \lambda & 0 \\
                                                      0 & \lambda^{-1} \\
                                                    \end{array}
                                                  \right),
                                                  \lambda=\pm 2,\lb{3.11}\\
N_1(\lambda,b)=\left(
 \begin{array}{cc}
  \lambda & b \\
   0& \lambda \\
   \end{array}
   \right),\lambda=\pm1,b=\pm1,0,\lb{3.12}\\
R(\theta)=\left(
            \begin{array}{cc}
              \cos{\theta} & -\sin{\theta} \\
              \sin{\theta} & \cos{\theta} \\
            \end{array}
          \right),\theta\in (0,\pi)\cup(\pi,2\pi),\lb{3.13}\\
N_2(\omega,B)=\left(
                \begin{array}{cc}
                  R(\theta) & B \\
                  0 &  R(\theta) \\
                \end{array}
              \right),\theta\in (0,\pi)\cup(\pi,2\pi),
\lb{3.14}\end{eqnarray}
where $B=\left(
           \begin{array}{cc}
             b_1 & b_2 \\
             b_3 & b_4 \\
           \end{array}
         \right)
$ with $b_i\in {\bf R}$ and $b_2\neq b_3$.

Splitting numbers possess the following properties:

{\bf Lemma 3.3.} (cf. \cite{Lon2}, Lemma 9.1.5 and List 9.1.12 of \cite{Lon4})
{\it For $M\in\Sp(2n)$, splitting numbers
$S^{\pm}_N(\omega)$ are constant for all $N \in \Omega^0(M)$. Moreover,
there hold
\begin{eqnarray}
S_M^{+}(\om) &=& S_M^{-}(\bar{\om}),~\forall~\omega
\in {\bf U}.\nn\\
S_M^{\pm}(\om) &=& 0, \qquad {\it if}\;\;\om\not\in\sg(M).  \nn\\
S_{N_1(1,a)}^+(1) &=& \left\{\matrix{1, &\quad {\rm if}\;\; a\ge 0, \cr
0, &\quad {\rm if}\;\; a< 0. \cr}\right. \nn\\S_{N_1(-1,a)}^+(-1) &=& \left\{\matrix{1, &\quad {\rm if}\;\; a\leq 0, \cr
0, &\quad {\rm if}\;\; a> 0. \cr}\right. \nn\\
(S_{R(\theta)}^+(e^{\sqrt{-1}\theta}), S_{R(\theta)}^-(e^{\sqrt{-1}\theta}))&=&(0,1), \theta \in (0, \pi) \cup (\pi, 2\pi) \nn\eea

For any $M_i\in\Sp(2n_i)$ with $i=0$ and $1$, there holds }
\bea S^{\pm}_{M_0\dm M_1}(\om) = S^{\pm}_{M_0}(\om) + S^{\pm}_{M_1}(\om),
    \qquad \forall\;\om\in\U. \nn\eea

{\bf Lemma 3.4.}(cf. Lemma 2.5 and Lemma 2.6 of \cite{LT2}) {\it For any $M\in Sp(2n)$ and $\omega\in {\bf U}$, the
splitting numbers $_PS^{\pm}_M (\omega)$ are well defined and
satisfy the following properties.

(i) $_PS^{\pm}_M (\omega)=S^{\pm}_{P^{-1}M}(\omega)-S^{\pm}_{P^{-1}}(\omega)$.

(ii) $_PS^{+}_M (\omega)= {}_PS^{-}_M (\bar{\omega})$.

(iii) $_PS^{\pm}_M (\omega)= {}_PS^{\pm}_N (\omega)$ if $P^{-1}N\in
\Omega^0(P^{-1}M)$.

(iv) $_PS^{\pm}_{M_1\diamond
M_2}(\omega)= {}_{P_1}S^{\pm}_{M_1}(\omega)+ {}_{P_2}S^{\pm}_{M_2}(\omega)$
for $M_j, P_j\in Sp(2n_j)$ with $n_j\in\{1,\cdots,n\}$ satisfying
$P=P_1\diamond P_2$ and $n_1+n_2=n$.

(v) $_PS^{\pm}_M (\omega)=0$ if $\omega\notin \sigma(P^{-1}M)\cup\sigma(P^{-1})$.}

{\bf Lemma 3.5.} {\it Assume $A(t)\in GL(\R^{2n})$ is positive definite for $t\in[0,\tau]$, let $\gamma\equiv\gamma_A\in \mathcal {P}_{\tau} (2n)$
be the fundamental solution of the linearized Hamiltonian system $\dot{y}(t)=JA(t)y(t)$. Then we have}
\bea i_\omega^P(\gamma)=\nu_\omega(P^{-1})+\sum_{0<s<\tau}\nu_\omega^P(\gamma(s)), \forall \omega\in\U.\nn\eea
{\bf Proof.} Similar to the proof of Proposition 15.1.3 of \cite{Lon4} and (4.4) of \cite{DoL1}, for any $\gamma\in \mathcal {P}_{\tau} (2n)$, we have
\bea i_\omega(\gamma|_{[0,s^+]})-i_\omega(\gamma|_{[0,s^-]})=\nu_\omega(\gamma(s)),\lb{jbds}\eea
if $A(t)=-J\dot{\gamma}(t)\gamma(t)^{-1}$ is positive definite on $[s-\epsilon, s+\epsilon]$ with $\epsilon>0$ is small.
Direct calculation give\bea B(t):=-J\frac{d}{dt}(\gamma_A(t)P^{-1})(\gamma_A(t)P^{-1})^{-1}&=&-J\dot{\gamma}_A(t)P^{-1}P\gamma_A(t)^{-1}\nn\\
&=&-J\dot{\gamma}_A(t)\gamma_A(t)^{-1}=A(t).\nn\eea
From (\ref{jbds}) one has
\bea i_\omega((P^{-1}\gamma_A)*\xi)&=&i_\omega(P(P^{-1}\gamma_A)P^{-1}*P\xi P^{-1})=i_\omega(\gamma_AP^{-1}*P\xi P^{-1})\nn\\
&=&i_\omega(\gamma_1P^{-1}*P\xi P^{-1})+\sum_{0<s<\tau}\nu_\omega(\gamma_A(s)P^{-1}),\lb{jbds2}\eea
where $\gamma_1(t):=\gamma(\epsilon t)$ for $t\in [0, \tau]$ and $\epsilon>0$ is small enough such that $\nu_\omega(\gamma(\epsilon t)P^{-1})=0$
for $t\in (0, \tau]$. If $\nu_\omega(\gamma(0)P^{-1})=\nu_\omega(P^{-1})=0$, then $i_\omega(\gamma_1P^{-1}*P\xi P^{-1})=i_\omega(P\xi P^{-1})=i_\omega(\xi)$.
If $\nu_\omega(P^{-1})\neq0$, then $\nu_\omega(\zeta_1(t))=0$ for $t\in(0, \tau]$, where $\zeta_1(t)=P^{-1}R(-\epsilon t)^{\diamond n}$ for
$t\in[0, \tau]$. Let $\zeta_2=(\gamma_1P^{-1})*\zeta_1^{-1}$, then $-J\dot{\zeta_2}(t)\zeta_2(t)^{-1}$ is positive definite
for $t\in[0, \tau]$, where for any path $\beta: [a, b]\rightarrow Sp(2n)$ we define $\beta^{-1}(t)=\beta(a+b-t)$ for $t\in[a, b]$, so
\bea i_\omega(\gamma_1P^{-1}*P\xi P^{-1})&=&\nu_\omega(P^{-1})+i_\omega(\zeta_1*P\xi P^{-1})\nn\\&=&\nu_\omega(P^{-1})+i_\omega(P\xi P^{-1})
=\nu_\omega(P^{-1})+i_\omega(\xi),\lb{jbds3}\eea
Combining (\ref{jbds2}) with (\ref{jbds3}), we have
\bea i_\omega^P(\gamma)&=&i_\omega((P^{-1}\gamma_A)*\xi)-i_\omega(\xi)\nn\\&=&\nu_\omega(P^{-1})+\sum_{0<s<\tau}\nu_\omega(\gamma_A(s)P^{-1})
\nn\\&=&\nu_\omega(P^{-1})+\sum_{0<s<\tau}\nu_\omega^P(\gamma(s)),\nn\eea
which completes the proof.\hfill\hb

Now, combining Lemma 3.5 with Propositions 2.4-2.5, we obtain a relationship between Ekeland P-index with Maslov P-index:

{\bf Theorem 3.6.} {\it Under the same assumption of Lemma 3.5, we have}
\bea i_1^P(\gamma)=\nu_1(P^{-1})+i_P^E(A).\nn\eea

\setcounter{equation}{0}
\section{Proof of the main result}

Define two function spaces $W_p$ and $L^2$ by:
\begin{eqnarray*}
W_p&=&\{x\in W^{1,2}([0,1],\mathbf{R}^{2n})\;|\;x(1)=Px(0)\}, \\
 L^2 &=& L^2((0,1),\mathbf{R}^{2n}) .
\end{eqnarray*}
And define a map $\Lambda:W_p\subset L^2\rightarrow L^2$ by $(\Lambda x)(t)=\dot{x}(t)$. We know $\Lambda$ is invertible, and for any $u\in L^2$ we obtain
\[
(\Lambda^{-1}u)(t)=x(t),
\]
\[
x(t)=\int_0^t u(\tau)\,d\tau+(P-I)^{-1}\int_0^1 u(\tau)\,d\tau.
\]
Note that $\Lambda^{-1}=\Pi_1$ and $x(1)=Px(0)$, where $\Pi_1$ is defined as in Section 2. So by Lemma 2.1,
we have $J\Lambda^{-1}:L^2\rightarrow L^2$ is self-adjoint and compact.

Consider the dual functional
\[
 \psi(u)=\int_0^1 \left[\frac{1}{2}(Ju,\Lambda^{-1}u)+H^*(-Ju)\right]\,dt
\]
on $L^\beta=L^2\cap L^\beta ((0,1);\mathbf{R}^{2n})$, with $\alpha^{-1}+\beta^{-1}=1$,
where $H^*(x^*)=\sup_{x\in \mathbf{R}^{2n}}\{(x,x^*)-H(x)\}$ is the Legendre transform of $H$ by Definition
\uppercase\expandafter{\romannumeral2.1.7} in \cite{Eke2}. As we know the global minimum of $\psi$ on $L^\beta$ is reached.
Next we will prove that $\bar{u}=\dot{\bar{x}}$ minimizes $\psi$
 \begin{equation}\label{80}
   \psi(\bar{u})=\inf\psi,
 \end{equation}
 where $\bar{x}$ is a solution of the boundary value problem:
\begin{equation}\label{70}
\left\{
\begin{array}{l}
\dot{x}(t)=JH'(x),\quad t\in (0,1)\\
 x(1)=Px(0).
\end{array}
\right.
\end{equation}

In fact, for $\forall u\in L^\beta$, we have $\psi'(u)\in L^\alpha$ is a linear functional on $L^\beta$:
\begin{eqnarray}
\psi'(u)(v)&=&\int_0^1 [(-J\Lambda^{-1}u,v)+(H^*(-Ju),-Jv)]\;dt,\qquad \forall u,v\in L^\beta. \nn\\
&=&\int_0^1 (-J\Lambda^{-1}u+JH^*(-Ju),v)\;dt. \label{90}
\end{eqnarray}
Because of $\psi'(\bar{u})=0$, we obtain:
\begin{equation}\label{81}
-J\Lambda^{-1}\bar{u}+JH^*(-J\bar{u})=0\quad \Rightarrow \quad H^*(-J\bar{u})=\Lambda^{-1}\bar{u}
\end{equation}
By the Legendre reciprocity formula of Proposition \uppercase\expandafter{\romannumeral2.1.15} and
Proposition \uppercase\expandafter{\romannumeral2.2.10} in \cite{Eke2}, we have
\[
  -J\bar{u}=H'(\Lambda^{-1}\bar{u}).
\]
Let $\bar{x}=\Lambda^{-1}\bar{u}$, then $\bar{u}=\dot{\bar{x}}$. It is clearly that $\bar{x}$ is a solution of the problem (\ref{70}).
Thus, we finish the proof.

From (\ref{90}), we have
\[
(\psi''(\bar{u})v,v)=\int_0^1 [(-J\Lambda^{-1}v,v)+(H^{*''}(-J(\bar{u}))Jv,Jv)]\; dt \quad \forall v\in L^\beta.
\]

By Proposition \uppercase\expandafter{\romannumeral2.2.10} in \cite{Eke2} and formula (\ref{81}), we get
\[
  I_{2n}=H^{*''}(-J\bar{u})H^{''}(\bar{x}) \quad \Rightarrow \quad H^{*''}(-J\bar{u})=(H^{''}(\bar{x}))^{-1}.
\]
Therefore, $\psi^{''}(\bar{u})$ can be defined on $L^2$. As $\bar{u}$ is a  minimal point of $\psi(u)$ by (4.1), we get
that Morse index of $\psi^{''}(\bar{u})$ defined on $L^2$ is zero, which means that
\bea i_P^E(A)=0, \lb{zeroindex}\eea
where $A=H^{\prime\prime}(\bar{x}(t))$
for $t\in[0,1]$, $i_P^E(A)$ is defined as in Section 2.

As $\bar{x}$ is defined on [0,1], we need to extend the domain to $[0,k]$. Let
\[
\bar{x}(t)=P^i\bar{x}(t-i), \qquad \forall t\in[i,i+1],~i=0,1,\cdots, k-1.
\]
By formula (\ref{81}), we obtain \[\lim_{\varepsilon\rightarrow 0^+} \bar{x}(1+\varepsilon)=\lim_{\varepsilon\rightarrow 0^+}
P\bar{x}(0+\varepsilon)=P\bar{x}(0)=\bar{x}(1).\]

So $\bar{x}\in C([0,2],\mathbf{R}^{2n})$, and $\bar{x}(2)=P\bar{x}(1)=P^2\bar{x}(0)$. By definition, we have
\[H(Py)=H(y), \quad \forall y\in \mathbf{R}^{2n}.\]
Thus we get
\begin{equation}\label{4.7}
P^{T}H'(Py)=H'(y) \quad \Rightarrow \quad PH'(y)=H'(Py).
\end{equation}
\begin{equation}\label{4.8}
P^{T}H''(Py)P=H''(y) \quad \Rightarrow \quad H''(Py)P=PH''(y), \quad \forall y \in \mathbf{R}^{2n}.
\end{equation}

We have from (\ref{4.7})
\begin{eqnarray*}
\lim_{\epsilon\rightarrow 0^+}\dot{\bar{x}}(1+\epsilon)&=&\lim_{\epsilon\rightarrow
0^+}P\dot{\bar{x}}(0+\epsilon)=\lim_{\epsilon\rightarrow 0^+}PJH'(\bar{x}(\epsilon))\\
&=&PJH'(\bar{x}(0))=JPH'(\bar{x}(0))=JH'(P\bar{x}(0))\\
&=&JH'(\bar{x}(1))=\lim_{\epsilon\rightarrow 0^+}JH'(\bar{x}(1-\epsilon))\\
&=&\lim_{\epsilon\rightarrow 0^+}\dot{\bar{x}}(1-\epsilon).
\end{eqnarray*}

Hence $\bar{x}\in C^1([0,2],\mathbf{R}^{2n})$, and $\bar{x}(t)$ satisfies
\[
\left\{
\begin{array}{l}
\dot{\bar{x}}(t)=JH'(\bar{x}),\quad \forall t\in (0,1)\\
\bar{x}(2)=P^2\bar{x}(0).
\end{array}
\right.
\]
By induction, we can finally get that $\bar{x}\in C^1([0,k],\mathbf{R}^{2n})$ and $\bar{x}(t)$ satisfies
\[
\left\{
\begin{array}{l}
\dot{\bar{x}}(t)=JH'(\bar{x}),\quad \forall t\in (0,k)\\
\bar{x}(k)=\bar{x}(0).
\end{array}
\right.
\]
Let $\gamma=\gamma_{\bar{x}}(t)$ be the fundamental solution of (\ref{1.3}) with $A(t)=H''(\bar{x}(t))$ for $t\in [0,k]$ satisfying $\gamma(0)=I_{2n}$.
By (\ref{4.8}), we get
\[
H''(\bar{x}(t+1))P=H''(P\bar{x}(t))P=PH''(\bar{x}(t)).
\]
Direct calculations give
\begin{eqnarray*}
\frac{d}{dt}\left(P\gamma(t)P^{-1}\gamma(1)\right)&=&P\dot{\gamma}(t)P^{-1}\gamma(1)\\
&=&PJH''(\bar{x}(t))\gamma(t)P^{-1}\gamma(1)\\
&=&JPH''(\bar{x}(t))\gamma(t)P^{-1}\gamma(1)\\
&=&JH''(\bar{x}(t+1))P\gamma(t)P^{-1}\gamma(1).
\end{eqnarray*}
Since $P\gamma(t)P^{-1}\gamma(1)\,|_{t=0}=\gamma(1)=\gamma(t+1)\,|_{t=0}$ and the fundamental solution of (\ref{1.3}) is unique, $\gamma$ satisfies
\[
\gamma(t+1)=P\gamma(t)P^{-1}\gamma(1), \qquad \forall t\in [0,k-1].
\]
Specially
\bea
\gamma(k)=P\gamma(k-1)P^{-1}\gamma(1)=\cdots=P^{k-1}\gamma(1)\left(P^{-1}\gamma(1)\right)^{k-1}=\left(P^{-1}\gamma(1)\right)^{k}.\lb{ddjz}
\eea
Combining (\ref{zeroindex}) with Theorem 3.6, note that $\ker(P-I_{2n})=0$,
we have \bea i_1^P(\gamma|_{[0,1]})=0,\lb{zeroindex2}\eea
By Lemma 3.5, there also holds \bea i_\omega^P(\gamma|_{[0,1]})\geq \nu_\omega(P^{-1}).\lb{omigaindex}\eea

{\bf Proof of Theorem 1.1.}

Note that $P\Sigma=\Sigma$ is equivalent to $P^{-1}\Sigma=\Sigma$.  We suppose that
$P=R(\theta_1)\diamond R(\theta_2)\diamond\cdots \diamond R(\theta_n)$, where $0<\frac{\theta_i}{\pi}\leq 1$.
Denote the eigenvalues of $P^{-1}$ on the upper
semi-circle in $\U$ by $\omega_1, \omega_2,\cdots, \omega_q$ anticlockwise. Let $\omega_0=1$.
By the definitions of splitting numbers, we have
\bea i^P_{\omega_i}(\gamma|_{[0,1]})+{}_PS^{+}_M (\omega_i)=i^P_{\omega_{i+1}}(\gamma|_{[0,1]})+{}_PS^{-}_M (\omega_{i+1}),i=0,1,...,q-1,\eea
where $M=\gamma(1)$.
Note that by Lemma 3.4(i), it follows that\bea {}_PS^{\pm}_M (\omega)=S^{\pm}_{P^{-1}M}(\omega)-S^{\pm}_{P^{-1}}(\omega).\eea
By Lemma 3.3, we have \bea S^{-}_{P^{-1}}(\omega_i)=S^{+}_{P^{-1}}(\bar{\omega}_i)=0,\quad
S^{+}_{P^{-1}}(\omega_i)=S^{-}_{P^{-1}}(\bar{\omega}_i)=\nu_{P^{-1}}(\omega_i),\quad if \quad \omega_i\neq -1,1,\nn
\\S^{-}_{P^{-1}}(\omega_i)=\frac{1}{2}\nu_{P^{-1}}(\omega_i)=
S^{+}_{P^{-1}}(\omega_i),\quad if \quad \omega_i= -1, 1.\eea
Note that $\omega_i\neq -1$ for $1 \leq i< q$ and $\nu_{P^{-1}}(1)=0$. Inserting (4.12)-(4.13) into (4.11), we obtain
\bea i^P_{\omega_1}(\gamma|_{[0,1]})-i^P_{1}(\gamma|_{[0,1]})&=&
S^{+}_{P^{-1}M}(1)-S^{-}_{P^{-1}M}(\omega_{1}),\\i^P_{\omega_{i+1}}(\gamma|_{[0,1]})-i^P_{\omega_i}(\gamma|_{[0,1]})&=&
S^{+}_{P^{-1}M}(\omega_i)-S^{-}_{P^{-1}M}(\omega_{i+1})-\nu_{P^{-1}}(\omega_i), \nn\\&& if\quad 1\leq i<q, \omega_q\neq -1,\\
i^P_{\omega_q}(\gamma|_{[0,1]})-i^P_{\omega_{q-1}}(\gamma|_{[0,1]})&=&
S^{+}_{P^{-1}M}(\omega_{q-1})-S^{-}_{P^{-1}M}(\omega_{q})\nn\\&&-\nu_{P^{-1}}(\omega_{q-1})+\frac{1}{2}\nu_{P^{-1}}(-1), \quad if \quad\omega_q=-1. \eea
Noticing that $\nu_{P^{-1}}(-1)=0$ if $\omega_q\neq -1$, thus (4.15)-(4.16) imply \bea i^P_{\omega_q}(\gamma|_{[0,1]})-i^P_{\omega_{q-1}}(\gamma|_{[0,1]})=
S^{+}_{P^{-1}M}(\omega_{q-1})-S^{-}_{P^{-1}M}(\omega_{q})-\nu_{P^{-1}}(\omega_{q-1})+\frac{1}{2}\nu_{P^{-1}}(-1),\eea
whenever $\omega_q$ is equal to $-1$ or not. Combining (4.14)-(4.15) and (4.17), we obtain
 \bea i^P_{\omega_{q}}(\gamma|_{[0,1]})-i^P_{1}(\gamma|_{[0,1]})&=&\sum_{i=0}^{q-1}(i^P_{\omega_{i+1}}(\gamma|_{[0,1]})-i^P_{\omega_i}(\gamma|_{[0,1]}))\nn\\
&=&\sum_{i=0}^{q-1}S^{+}_{P^{-1}M}(\omega_i)-\sum_{i=1}^{q}S^{-}_{P^{-1}M}(\omega_{i})\nn\\&&-\sum_{i=1}^{q-1}\nu_{P^{-1}}(\omega_i)+\frac{1}{2}\nu_{P^{-1}}(-1),\eea
which together with (\ref{zeroindex2})-(\ref{omigaindex}) implies
\bea \sum_{i=0}^{q-1}S^{+}_{P^{-1}M}(\omega_i)\geq \left\{\matrix{
            \sum_{i=1}^{q-1}\nu_{P^{-1}}(\omega_i)+\frac{1}{2}\nu_{P^{-1}}(-1)=n, & \quad {\rm if}\;\omega_q=-1, \cr
            \sum_{i=1}^{q}\nu_{P^{-1}}(\omega_i)=n, & \quad {\rm if}\; \omega_q\neq-1. \cr}\right.
\eea
On the other hand, note that $S^{+}_{P^{-1}M}(\omega_0)=S^{+}_{P^{-1}M}(1)\leq\frac{1}{2}\nu_{P^{-1}M}(1)$
and $S^{+}_{P^{-1}M}(\omega_i)\leq \nu_{P^{-1}M}(\omega_i)$, we have
\bea \sum_{i=0}^{q-1}S^{+}_{P^{-1}M}(\omega_i)\leq \frac{1}{2}\nu_{P^{-1}M}(1)+\sum_{i=1}^{q-1}\nu_{P^{-1}M}(\omega_i)\leq \frac{1}{2}e(P^{-1}M)\leq n,\eea
Comparing (4.20) with (4.19), we obtain $e(P^{-1}\gamma(1))=e(P^{-1}M)=2n$, which together with (\ref{ddjz})
gives $e(\gamma(k))=2n$, thus $\bar{x}$ corresponds to an elliptic closed characteristic on $\Sigma$.\hfill\hb

{\bf Proof of Theorem 1.2.}

As in Section 1.7
of \cite{Eke2}(cf. also Proposition 2.13 of \cite{Liu1}), we have \begin{eqnarray}P^{-1}\gamma(1)=P_{\bar{x}}^{-1}(N_1(1,1)
\diamond
Q)P_{\bar{x}}
\end{eqnarray}
for some symplectic matrices $P_{\bar{x}}$ and $Q$(cf. also Lemmas 15.2.3 and 15.2.4 of \cite{Lon4}), where $N_1(1,1)=\left(
                    \begin{array}{cc}
                      1 & 1 \\
                      0 & 1 \\
                    \end{array}
                  \right)$. Now we proceed the proof by contradiction.
We suppose that $\bar{x}$ corresponds to a hyperbolic closed characteristic on $\Sigma$, then $e(P^{-1}\gamma(1))=2$
and $\sigma(Q)\cap \U=\emptyset$ by (4.21), where $\sigma(Q)$ denotes the spectrum of $Q$. Thus by Lemma 3.3, we obtain
\bea S^{-}_{P^{-1}\gamma(1)}(1)=0, \quad S^{+}_{P^{-1}\gamma(1)}(1)=1,\nn\\
S^{-}_{P^{-1}\gamma(1)}(\omega)=S^{+}_{P^{-1}\gamma(1)}(\omega)=0, \quad if \quad \omega\in \U\setminus\{1\}.\eea
Without loss of generality, we suppose $P=R(\theta_1)\diamond R(\theta_2)\diamond\cdots \diamond R(\theta_n)$,
where $0<\frac{\theta_i}{\pi}< 1$ for $1\leq i\leq a$, $\frac{\theta_i}{\pi}=1$ for $a+1\leq i\leq a+b$, $1<\frac{\theta_i}{\pi}< 2$ for $a+b+1\leq i\leq n$
$a,b\geq 0$. By the assumption that $^\#\{i\mid \theta_i\in(0, \pi]\}-^\#\{i\mid \theta_i\in(\pi, 2\pi)\}\geq 2$, we have \bea a+b-(n-a-b)\geq2.\eea
Denote the eigenvalues of $P^{-1}$ on the upper
semi-circle in $\U$ by $\omega_1, \omega_2,\cdots, \omega_q$ anticlockwise. Let $\omega_0=1$.
By the definitions of splitting numbers, we get
\bea i^P_{\omega_i}(\gamma|_{[0,1]})+{}_PS^{+}_M (\omega_i)=i^P_{\omega_{i+1}}(\gamma|_{[0,1]})+{}_PS^{-}_M (\omega_{i+1}),i=0,1,...,q-1,\eea
where $M=\gamma(1)$. Note that by Lemma 3.4(i), it follows that\bea {}_PS^{\pm}_M (\omega)=S^{\pm}_{P^{-1}M}(\omega)-S^{\pm}_{P^{-1}}(\omega).\eea
Combining (4.22) with (4.24)-(4.25), we obtain\bea i^P_{\omega_{1}}(\gamma|_{[0,1]})-i^P_{\omega_0}(\gamma|_{[0,1]})&=&1+S^{-}_{P^{-1}}(\omega_{1})\nn\\
i^P_{\omega_{i+1}}(\gamma|_{[0,1]})-i^P_{\omega_i}(\gamma|_{[0,1]})&=&
-S^{+}_{P^{-1}}(\omega_i)+S^{-}_{P^{-1}}(\omega_{i+1}), \quad 1\leq i<q.\eea
By Lemma 3.3, we have \bea S^{-}_{P^{-1}}(\omega)=S^{+}_{P^{-1}}(\bar{\omega})=0,\quad
S^{+}_{P^{-1}}(\omega)=S^{-}_{P^{-1}}(\bar{\omega})=\nu_{P^{-1}}(\omega),\quad if \quad \omega=e^{\sqrt{-1}\theta_i}, 1\leq i\leq a,\nn
\\S^{-}_{P^{-1}}(\omega)=\frac{1}{2}\nu_{P^{-1}}(\omega)=
S^{+}_{P^{-1}}(\omega),\quad if \quad \omega= -1,\eea
which together with (4.26) implies
\bea i^P_{\omega_{q}}(\gamma|_{[0,1]})-i^P_{\omega_0}(\gamma|_{[0,1]})&=&
1+\sum_{i=1}^{q}S^{-}_{P^{-1}}(\omega_i)-\sum_{i=1}^{q-1}S^{+}_{P^{-1}}(\omega_{i})\nn
\\&=&\left\{\matrix{
            1+(n-a-b+b)-a, & \quad {\rm if}\;\omega_q=-1, \cr
           1+(n-a-b)-(a-S^{+}_{P^{-1}}(\omega_{q})), & \quad {\rm if}\; \omega_q\neq-1. \cr}\right.\eea
On the other hand, it follows from (\ref{zeroindex2})-(\ref{omigaindex}) and (4.27) that
\bea i^P_{\omega_{q}}(\gamma|_{[0,1]})-i^P_{\omega_0}(\gamma|_{[0,1]})&\geq&\nu_{P^{-1}}(\omega_q)
\nn
\\&=&\left\{\matrix{
            2b, & \quad {\rm if}\;\omega_q=-1, \cr
           S^{+}_{P^{-1}}(\omega_{q}), & \quad {\rm if}\; \omega_q\neq-1. \cr}\right.\eea
By (4.23), we have $2b>1+(n-a-b+b)-a$ and $S^{+}_{P^{-1}}(\omega_{q})>1+(n-a-b)-(a-S^{+}_{P^{-1}}(\omega_{q}))$
when $\omega_q\neq-1$ which implies $b=0$. Comparing (4.29) with (4.28), we get a contradiction and complete the proof. \hfill\hb

\bibliographystyle{abbrv}

\end{document}